\documentclass[12pt]{article}
\usepackage{amsmath,amsthm,amssymb}
\usepackage{mathtools}
\newtheorem{theorem}{Theorem}[section] 
\newtheorem{lemma}{Lemma}[section]
\newtheorem{cor}{Corollary}[section]

\newtheorem{ass}{Assumption}[section]

\newtheorem{proposition}{Proposition}[section]

\numberwithin{equation}{section}

\newcommand{\D}{{\rm d}}
\newcommand{\dx}{\, \D x}

\newcommand{\rz}{\mathbb{R}}
\newcommand{\nz}{\mathbb{N}}

\newcommand{\eps}{\varepsilon}

\newcommand{\klauf}{\left(\begin{array}}
\newcommand{\klzu}{\end{array}\right)}

\title{Small weights in Caccioppoli's inequality and  
applications to Liouville-type theorems for non-standard problems}

\author{Michael Bildhauer \& Martin Fuchs}
\date{}

\newcommand{\reff}[1]{(\ref{#1})}

\usepackage{hyperref}

\hypersetup{
   pdfnewwindow=true,
   colorlinks=false,
   linkcolor=black,
   citecolor=green,
   filecolor=magenta,
   urlcolor=cyan,
}

\begin{document}

\parindent0em
\maketitle

\newcommand{\op}[1]{\operatorname{#1}}
\newcommand{\bv}{\op{BV}}
\newcommand{\mub}{\overline{\mu}}
\newcommand{\muhat}{\hat{\mu}}

\newcommand{\hypref}[2]{\hyperref[#2]{#1 \ref*{#2}}}
\newcommand{\hypreff}[1]{\hyperref[#1]{(\ref*{#1})}}

\newcommand{\ob}[1]{^{(#1)}}

\newcommand{\xh}{\Xi}
\newcommand{\oh}[1]{O\left(#1\right)}
\newcommand{\xn}{w_{1}}
\newcommand{\yn}{w_{2}}
\newcommand{\On}{\hat{\Omega}}
\newcommand{\uq}{\hat{u}}
\newcommand{\gameps}{\Gamma_{Q}}
\newcommand{\se}{{s_Q}}
\begin{abstract}
Using a variant of Caccioppoli's inequality involving small weights, i.e.~weights of the form $(1+|\nabla u|^2)^{-\alpha/2}$ 
for some $\alpha > 0$, we establish several
Liouville-type theorems under general non-standard growth conditions.{}
\end{abstract}

\vspace*{2ex}

\centerline{Dedicated to the 70th birthday of Gregory Seregin} 

\vspace*{2ex}

\parindent0ex
%***************************************************************************************************
%***************************************************************************************************
\section{Introduction}\label{intro}
%***************************************************************************************************
%***************************************************************************************************

%***** LABEL intro .... *****
 
 Throughout this manuscript we always suppose that $u$: $\rz^2 \to \rz$, $u \in C^2(\rz^2)$, is a solution of the nonlinear equation
\begin{equation}\label{intro 1}
\op{div} \big[\nabla f(\nabla u)\big] = 0 \quad\mbox{on $\rz^2$}\, .
\end{equation} 
 
 Our main goal is the discussion of Liouville-type results under rather general hypotheses on the convex density
 $f$: $\rz^2 \to \rz$ including non-standard growth conditions such as the case of linear growth or even
 allowing a certain degree of anisotropy in the superlinear case. For technical simplicity we restrict ourselves to the twodimensional  case.\\ 
 
 It is out of reach to give a complete overview on all the recent contributions on Liouville-type results.
 We refer to the beautiful survey of Farina \cite{Fa:2007_1} in the case of general elliptic problems
including a lot of historical references. We also refer to Seregin's discussion \cite{Se:2018_1} of Liouville-type theorems
in the case of the Navier-Stokes equations.\\

Contributions in the case of linear or anisotropic growth are quite rarely found. We just mention the papers
\cite{Am:2009_1}, \cite{AG:2015_1}, \cite{Du:2020_1}, \cite{BF:2021_1} and the references quoted therein.\\

Before going into details we fix our main assumption which will be supplemented with appropriate hypotheses
adapted to the applications of \mbox{Section \ref{linear}} and of Section \ref{pq}.

\begin{ass}\label{ass main 1}
The convex energy density $f$: $\rz^2 \to \rz$ is of class $C^2(\rz^2)$ and satisfies the non-uniform
 ellipticity condition
\begin{equation}\label{intro 2}
c_1 \big(1+|Z|^2\big)^{-\frac{\mu}{2}} |Y|^2 \leq D^2 f(Z)\big(Y,Y\big) \leq c_2 \big(1+|Z|^2\big)^{-\frac{\mub}{2}}|Y|^2
\end{equation}
for all $Z$, $Y\in \rz^2$ with exponents $\mu > 1$, $\mub \leq 1$ and constants $c_1$, $c_2 >0$.
\end{ass}

Condition \reff{intro 2} also serves as one main assumption in the recent paper \cite{BF:2021_1} on Liouville-type results
in two dimensions for functionals satisfying a linear growth condition. These results are restricted to the case $\mub =1$, for instance we have:

\begin{theorem}\label{intro theo 1} (Theorem 1.1, c), \cite{BF:2021_1}, the case $N=1$)
Let $u \in C^2(\rz^2)$ denote a solution of \reff{intro 1} with density $f$ such that for some real number $M>0$ we have 
\[
|\nabla f(Z)| \leq M \qquad\mbox{for all}\quad Z \in \rz^2
\]
and such \reff{intro 2} holds with the choice $\mub =1$.\\

If we have
\begin{equation}\label{intro 3}
\limsup_{|x| \to \infty} \frac{|u(x)|}{|x|} < \infty \, ,
\end{equation}
then $u$ is affine.
\end{theorem}

Let us give some further explanations concerning the hypotheses of Theorem \ref{intro theo 1}: from \reff{intro 2} and the boundedness
of $\nabla f$ it follows that $f$ is of linear growth, and this actually holds for arbitrary exponents $\mub \leq 1 < \mu$. We refer
to Lemma 2.1 in \cite{BF:2021_1}.\\

At the same time, if $f$ is of linear growth satisfying inequality \reff{intro 2}, then according to Lemma 1.1 of \cite{BF:2020_3} we necessarily
get the restriction $\mub \leq 1$, whereas the bound $\mu >1$ is an immediate consequence of the linear growth of $f$. To sum up,
Theorem \ref{intro theo 1} addresses the case of energy densities with linear growth, but just covers the
``limit case'' for which $\mub =1$. So the first question arises, whether the result of Theorem \ref{intro theo 1}
keeps valid, if we allow exponents $\mub < 1$.\\

Closely related is the following setting in the case of superlinear growth. Suppose that we have \reff{intro 2} with exponent
$\mub = 2-q$, $q >1$, on the right-hand side. Then we are interested in the anisotropic case, which means that we 
do not narrow our discussion by assuming $q$-power growth of $f$. We just 
impose the inequality
\[
a |Z|^s -b \leq f(Z) \qquad\mbox{for all}\quad Z \in \rz^2\, ,
\]
with exponent $1 < s \leq q$ and with constants $a$, $A >0$, $a$, $b\geq 0$ as lower bound for the density $f$ (see Section \ref{pq} for some further
comments on the assumptions). A Liouville-type result in this setting is established in Section \ref{pq}.
We emphasize that we are not aware of similar Liouville-type theorems w.r.t.~this general kind of anisotropic
hypotheses.\\

Let us fix the notation.\\

{\bf Notation.} We always abbreviate
\[
\Gamma := 1+|\nabla u|^2\, , \quad \gameps := 1 +|\nabla u - Q|^2 \qquad\mbox{for vectors} \quad Q\in \rz^2 \, .
\]
For fixed radii $r$, $R>0$ we consider open disks $B_r$ and define
\[
T_R := B_{2R} - \overline{B_{R}} \, ,\quad \hat{T}_{R}: = B_{5R/2} - \overline{B_{R/2}} \, ,
\]
where the center $x_0$ is not indicated. \\

Moreover, for $r >0$ we let
\[
(\nabla u) = (\nabla u)_r =-\hspace*{-2.5ex} \int_{B_r}\nabla u \dx \, .
\]

Then we have
\begin{theorem}\label{theo 1} Assume that  Assumption \ref{ass main 1} holds. Moreover, suppose
that there are real numbers $\gamma$, $\gamma_Q \geq 0$ such that
\begin{equation}\label{theo ass 1}
\gamma + \gamma_Q  <  \frac{1}{2}
\end{equation}
and that there exist $Q = Q(R) \in \rz^2$ such that
\begin{equation}\label{theo ass 2}
\liminf_{R \to \infty} \frac{1}{R^2} \Xi(R) :=
\liminf_{R \to \infty} \frac{1}{R^2} \int_{T_R}  \Gamma^{-\frac{\gamma + \mub}{2}}\gameps^{-\frac{\gamma_Q}{2}} |\nabla u-Q|^{2} \dx
< \infty \, .
\end{equation}
Then $u$ is an affine function.  
\end{theorem}

We note that Theorem \ref{theo 1} just relies on condition \reff{intro 2} and it does not matter, whether 
the energy density $f$ is of linear growth or even shows a completely anisotropic behaviour.  
Moreover, the conclusion of the theorem is independent of the value of the exponent $\mu$.\\

Let us shortly comment on the main idea for the proof of Theorem \ref{theo 1} recalling the approach towards
Theorem \ref{intro theo 1}. This theorem is proved by first using a Caccioppoli-type inequality
for the differentiated Euler equation. Since we have $\mub =1$ as a hypothesis of Theorem \ref{intro theo 1}, the right-hand side of this inequality 
can be measured in terms of the quantity 
$\nabla f(\nabla u) \cdot \nabla u$ which in turn occurs on the left-hand side of the weak form of \reff{intro 1} applied to a suitable
testfunction. \\

In the case $\mub < 1$ a serious gap arises which cannot be closed by obvious arguments. Here, as the main new feature, we introduce small weights
in Caccioppoli's inequality such that both sides again fit together. The same arguments bridge the gap in the superlinear
anisotropic case.\\

Theorem \ref{theo 1} immediately gives the following elementary corollary.

\begin{cor}\label{cor 1}
Suppose that we have Assumption \ref{ass main 1}. Then $u$ is an 
affine function if one of the following conditions hold.
\begin{enumerate}
\item $\displaystyle \sup_{x\in \rz^2} |\nabla u(x)| < \infty \, .$

\item There exists some $\eps >0$ such that
\[
 -\hspace*{-2.5ex}\int_{T_{R}} \Gamma^{\frac{1}{2}\left(\frac{3}{2}-\mub\right)+\eps}\dx < c
\] 
 with a constant $c>0$ not depending on $R$.\\
\end{enumerate}
\end{cor}

Before formulating more refined corollaries, we establish our Caccioppoli-type estimate in the next section as the main tool for 
proving Theorem \ref{theo 1} in Section \ref{proof}.\\

Section \ref{linear} is devoted to applications in the linear growth setting
while Section \ref{pq} concentrates on the main corollary in the superlinear case.\\

We finally note that the generality of Theorem \ref{theo 1} may be used to discuss a series of other applications
which is left to the particular interest of the reader.

%**************************************************************************************************
%**************************************************************************************************
\section{A Caccioppoli-type inequality}\label{cacc}
%**************************************************************************************************
%**************************************************************************************************

We start with a Caccioppoli-type inequality weighted with negative powers of $\Gamma$ and $\gameps$.

%***** LABEL cacc .... *******************************************************************************
\begin{lemma}\label{cacc lem 1}
Given Assumption \ref{ass main 1} we fix $Q\in \rz^2$, consider real numbers $s_Q> -1/4$, $s_1>-1/4$ and let
\[
c_Q := \left\{\begin{array}{ccl} 4 |\se|&\mbox{if}& \se\leq 0\, ,\\[1ex]
0&\mbox{if}& \se >0 \, ,\end{array}\right. \qquad
c_1 := \left\{\begin{array}{ccl} 4 |s_1|&\mbox{if}& s_1\leq 0\, ,\\[1ex]
0&\mbox{if}& s_1 >0 \, ,\end{array}\right. \qquad
\]
and suppose that  $\eta\in C^{\infty}_{0}\big(\rz^2\big)$, $0 \leq \eta \leq 1$. \\

If $c_Q + c_1 < 1$, then we have  (summation w.r.t.~$\alpha =1$, $2$)
\begin{eqnarray}\label{cacc 1}
\lefteqn{\big[1-c_Q-c_1\big] \int_{\rz^2} D^2 f (\nabla u) \big(\nabla \partial_\alpha u, \nabla \partial_\alpha u)
\gameps^{\se} \Gamma^{s_1}\eta^2 \dx}\nonumber\\
&\leq & c \Bigg[ \int_{\op{spt}\nabla \eta} D^2 f (\nabla u) \big(\nabla \partial_\alpha u, \nabla \partial_\alpha u)
\gameps^{\se} \Gamma^{s_1}\eta^2 \dx\Bigg]^{\frac{1}{2}}
\nonumber\\
&& \cdot \Bigg[ \int_{\op{spt}\nabla \eta} D^2f(\nabla u)\big(\nabla \eta, \nabla \eta\big) \big|\nabla u- Q\big|^2
\gameps^{\se}\Gamma^{s_1} \dx\Bigg]^{\frac{1}{2}}\, ,
\end{eqnarray}
where the constant $c$ is not sdepending on $\eta$.
In particular we have
\begin{eqnarray}\label{cacc 2}
\lefteqn{\int_{\rz^2} D^2 f (\nabla u) \big(\nabla \partial_\alpha u, \nabla \partial_\alpha u)
\gameps^{\se}\Gamma^{s_1} \eta^2\dx}\nonumber\\
&\leq& c  \int_{\op{spt}\nabla \eta} D^2f(\nabla u)\big(\nabla \eta, \nabla \eta\big) 
\big|\partial_\alpha u - Q_\alpha \big|^2 \gameps^{\se}\Gamma^{s_1} \dx \, .
\end{eqnarray}
\end{lemma} 

\vspace*{2ex}
\emph{Proof.} We first consider the case that $-1/4 < s_Q \leq 0$. For $\alpha =1$, $2$ and
for all $\psi \in C^{\infty}_{0}(\rz^2)$ equation \reff{intro 1} yields
\begin{equation}\label{cacc 3}
0 = \int_{\rz^2} D^2f(\nabla u) \big(\nabla \partial_\alpha u, \nabla \psi\big) \dx \, .
\end{equation}
Inserting $\psi: = \eta^2 \big(\partial_\alpha u- Q_\alpha\big) \gameps^\se \Gamma^{s_1}$ in \reff{cacc 3} 
we obtain for $\alpha =1$, $2$ and any testfunction $\eta$
\begin{eqnarray}\label{cacc 4}
\lefteqn{\int_{\rz^2} D^2f(\nabla u) \big(\nabla \partial_\alpha u ,\nabla \partial_\alpha u\big) 
\gameps^{\se}\Gamma^{s_1} \eta^2 \dx}\nonumber\\
&=& - \int_{\rz^2} D^2f(\nabla u)\big(\nabla \partial_\alpha u, \big(\partial_\alpha  u -Q_\alpha\big)
\nabla \gameps^{\se}\big)\Gamma^{s_1}\eta^2 \dx\nonumber\\
&& - \int_{\rz^2} D^2f(\nabla u)\big(\nabla \partial_\alpha u, \big(\partial_\alpha  u -Q_\alpha\big)
\nabla \Gamma^{s_1}\big)\gameps^{\se}\eta^2 \dx\nonumber\\
&& - 2 \int_{\rz^2} D^2 f(\nabla u) \big(\nabla \partial_\alpha u , \nabla \eta\big)
\big(\partial_\alpha u - Q_\alpha\big)  \gameps^{\se}\Gamma^{s_1}\eta \dx \, .
\end{eqnarray}

We denote the bilinear form $D^2f(\cdot,\cdot)$ by $\langle \cdot,\cdot\rangle$ and discuss the left-hand side of  \reff{cacc 4}
by observing that 
\begin{eqnarray}\label{cacc 5}
\lefteqn{\sum_{\alpha =1}^2 \Big\langle \partial_\alpha \nabla u, \partial_\alpha \nabla u\Big\rangle \gameps^\se}\nonumber\\
&\geq & \sum_{\alpha =1}^2 \Big\langle \partial_\alpha \nabla u, \partial_\alpha \nabla u\Big\rangle
\Big[\big(\partial_1 u - Q_1\big)^2 + \big(\partial_2 u -Q_2\big)^2\Big]   \gameps^{\se-1}\nonumber \\
& \geq & \sum_{\alpha =1}^2 \Big\langle \big(\partial_{\gamma} u -Q_\alpha\big) \partial_\alpha \nabla u , 
\big(\partial_\alpha u - Q_\alpha\big)  \partial_\alpha \nabla u \Big\rangle \gameps^{\se-1} \, .
\end{eqnarray}

Moreover, for the first integral on the right-hand side of \reff{cacc 4} we write
\begin{eqnarray}\label{cacc 6}
\lefteqn{\sum_{\alpha=1}^2\Big\langle \partial_\alpha \nabla u , \big(\partial_\alpha u -Q_\alpha\big)\nabla  \gameps^{\se}\Big\rangle}\nonumber\\
&=& \se \sum_{\alpha=1}^2 \Big\langle \big(\partial_\alpha u -Q_\alpha\big) \partial_\alpha \nabla u , 
\nabla \sum_{i=1}^2 \big(\partial_i u - Q_i\big)^2\Big\rangle 
 \gameps^{\se-1}\nonumber\\ 
 &=& 2\se \sum_{\alpha =1}^2 \Big\langle \big(\partial_\alpha u -Q_\alpha\big) \partial_\alpha \nabla u,  
 \big(\partial_\alpha u - Q_\alpha \big)\partial_\alpha \nabla u\Big\rangle
 \gameps^{\se-1}\nonumber\\
 && + 4\se \Big\langle\big( \partial_1 u -Q_1\big) \partial_1 \nabla u,  \big(\partial_2 u - Q_2\big)  \partial_2 \nabla u\Big\rangle
 \gameps^{\se-1} \, .
\end{eqnarray}

On account of
\begin{eqnarray*}
\lefteqn{\Big|\Big\langle \big(\partial_1 u -Q_1) \partial_1 \nabla u,  \big(\partial_2 u - Q_2\big) \partial_2 \nabla u\Big\rangle\Big|}\\
&\leq & \frac{1}{2} \sum_{\alpha =1}^2 \Big\langle \big(\partial_\alpha u -Q_\alpha\big) \partial_\alpha\nabla u , 
\big(\partial_\alpha u -Q_\alpha\big) \partial_\alpha \nabla u \Big\rangle
\end{eqnarray*}

we obtain from \reff{cacc 6}
\begin{eqnarray}\label{cacc 7}
\lefteqn{\Bigg|\sum_{\alpha=1}^2\Big\langle \partial_\alpha \nabla u , \big(\partial_\alpha u -Q_\alpha\big) 
\nabla  \gameps^{\se}\Big\rangle\Bigg| }
\nonumber\\
&& \leq 4 |\se| \sum_{\alpha =1}^2 \Big\langle \big(\partial_\alpha u -Q_\alpha\big) \partial_\alpha \nabla u, 
\big( \partial_\alpha u -Q_\alpha \big) \partial_\alpha \nabla u\Big\rangle
 \gameps^{\se-1} \, .
\end{eqnarray}

Combining \reff{cacc 7} and \reff{cacc 5} we get (where from now on we take the sum w.r.t.~$\alpha =1$, $2$)
\begin{eqnarray}\label{cacc 8}
\lefteqn{- \int_{\rz^2} D^2f(\nabla u)\big( \partial_\alpha \nabla u , \big(\partial_\alpha u-Q_\alpha\big)  
\nabla  \gameps^\se\big)\Gamma^{s_1}\eta^2 \dx} \nonumber\\
&&\leq 4 |\se| \int_{\rz^2} D^2f(\nabla u) \big(\partial_\alpha \nabla u, \partial \gamma \nabla u\big) \gameps^\se 
\Gamma^{s_1} \eta^2 \dx\, .
\end{eqnarray}

In the case $0 < \se$ we just use the positive sign of
\begin{equation}\label{cacc 9}
\Big\langle \partial_\alpha \nabla u , \big(\partial_\alpha u-Q_\alpha\big),  \nabla  \gameps^\se\Big\rangle
=\frac{1}{2}  \Big\langle \nabla \big(\partial_\alpha u -Q_\alpha\big)^2  , \nabla  \gameps^\se\Big\rangle \, .
\end{equation}

Having established \reff{cacc 8} and \reff{cacc 9} we recall $c_Q: = 4|\se|$ if $-1/4 < \se \leq 0$ and $c_Q = 0$ if $\se >0$. Then we summarize
\reff{cacc 8} and \reff{cacc 9} by writing
\begin{eqnarray}\label{cacc 10}
\lefteqn{- \int_{\rz^2} D^2f(\nabla u)\big( \partial_\alpha \nabla u , \big(\partial_\alpha u-Q_\alpha\big)  
\nabla  \gameps^\se\big)\Gamma^{s_1}\eta^2 \dx} \nonumber\\
&&\leq c_Q \int_{\rz^2} D^2f(\nabla u) \big(\partial_\alpha \nabla u, \partial \alpha \nabla u\big) \gameps^\se 
\Gamma^{s_1}\eta^2 \dx\, .
\end{eqnarray}

In the same way we recall $c_1 : = 4|s_1|$ if $-1/4 < s_1 \leq 0$ and $c_1 = 0$ if $s_1 >0$. With exactly the same reasoning as above
we additionally obtain
\begin{eqnarray}\label{cacc 11}
\lefteqn{- \int_{\rz^2} D^2f(\nabla u)\big( \partial_\alpha \nabla u , \big(\partial_\alpha u-Q_\alpha\big)  
\nabla  \Gamma^{s_1}\big)\gameps^{\se} \eta^2\dx} \nonumber\\
&&\leq c_1 \int_{\rz^2} D^2f(\nabla u) \big(\partial_\alpha \nabla u, \partial \alpha \nabla u\big) \gameps^\se 
\Gamma^{s_1} \eta^2\dx\, .
\end{eqnarray}

Returning to \reff{cacc 4} and using \reff{cacc 10} and \reff{cacc 11} we get
\begin{eqnarray}\label{cacc 12}
\lefteqn{\big[1-c_\eps-c_1\big]  \int_{\rz^n} D^2 f(\nabla u) \big(\nabla \partial_\alpha u , \nabla \partial_\alpha u\big) \gameps^s \eta^2 \dx}\nonumber\\
&&  \leq -2 \int_{\rz^n} D^2 f(\nabla u) \big(\eta \nabla \partial_\alpha u, \big(\partial_\alpha u -Q_\alpha\big)  \nabla \eta\big) 
\gameps^{\se} \Gamma^{s_1} \dx \, .
\end{eqnarray}

On the right-hand side of \reff{cacc 12} we observe that the integration is performed w.r.t.~the domain $\op{spt}\nabla \eta$ and 
an application of the
Cauchy-Schwarz inequality completes the proof of Lemma \ref{cacc lem 1}. \qed

%***************************************************************************************************
%***************************************************************************************************
\section{Proof of Theorem \ref{theo 1}}\label{proof}
%***************************************************************************************************
%***************************************************************************************************

%***** LABEL proof .... ******************************************************************************

For the proof of Theorem \ref{theo 1} we fix a disk $B_r \subset \rz^2$. 
We apply the Sobolev-Poincar\'{e} inequality to the solution $u \in C^2(\rz^2)$ under consideration and get the inequality
\begin{eqnarray}\label{proof 1}
\int_{B_r} \big| \nabla u - (\nabla u)\big|^2 \dx &\leq& c \Bigg[ \int_{B_r} \big| \nabla^2 u\big|\dx \Bigg]^2\nonumber \\[2ex]
& \leq & c\Bigg[ \int_{B_r}\Gamma^{-\frac{\mu}{4}} |\nabla^2 u| \Gamma^{\frac{\mu}{4}}\Gamma^{-\frac{\gamma}{4}}
\gameps^{-\frac{\gamma_Q}{4}}\Gamma^{\frac{\gamma}{4}}\gameps^{\frac{\gamma_Q}{4}}\dx\Bigg]^2 \nonumber \\[2ex]
&\leq & c \Bigg[\int_{B_r} \Gamma^{-\frac{\mu}{2}} |\nabla^2 u|^2 \Gamma^{-\frac{\gamma}{2}} \gameps^{-\frac{\gamma_Q}{2}}\dx\Bigg]
\, \Bigg[\int_{B_r} \Gamma^{\frac{\mu+\gamma}{2}}\gameps^{\frac{\gamma_Q}{2}} \dx \Bigg]\nonumber \\[2ex]
&\leq & c \int_{B_r} \Gamma^{-\frac{\mu}{2}} |\nabla^2 u|^2 \Gamma^{-\frac{\gamma}{2}}\gameps^{-\frac{\gamma_Q}{2}}\dx \, ,
\end{eqnarray}
where we used the fact that $|\nabla u|$ is bounded on the fixed ball $B_r$, however $c$ may depend on the radius $r$.\\

Now we choose $R \gg r$ and let $\eta \in C^\infty_{0}\big(B_{2R}\big)$, $0 \leq \eta \leq 1$, such that $\eta \equiv 1$ on $B_R$,
$|\nabla \eta| \leq c /R$. Then \reff{cacc 1} of Lemma \ref{cacc lem 1} gives recalling \reff{theo ass 1}
\begin{eqnarray}\label{proof 2}
\lefteqn{ \int_{B_r} \Gamma^{-\frac{\mu}{2}} |\nabla^2 u|^2 \Gamma^{-\frac{\gamma}{2}}\gameps^{-\frac{\gamma_Q}{2}}\dx}
\nonumber\\[1ex]
 &\leq &  \int_{B_{2 R}} \Gamma^{-\frac{\mu}{2}} |\nabla^2 u|^2 \Gamma^{-\frac{\gamma}{2}} 
 \gameps^{-\frac{\gamma_Q}{2}}\eta^2 \dx \nonumber\\[1ex]
 &\leq & c\Bigg[\int_{\op{spt}\eta} D^2f \big(\nabla \partial_\alpha u, \nabla \partial_\alpha u\big) \Gamma^{-\frac{\gamma}{2}}
\gameps^{-\frac{\gamma_Q}{2}} \eta^2\dx\Bigg]^{\frac{1}{2}}\nonumber\\[1ex]
 && \cdot\Bigg[ \frac{c}{R^2} \int_{T_R} \Gamma^{-\frac{\gamma +\mub}{2}} 
\gameps^{-\frac{\gamma_Q}{2}} \big|\partial_\alpha u - Q_\alpha\big|^2 \dx\Bigg]^{\frac{1}{2}}\nonumber\\[1ex]
& =:&c  I(R) \cdot \Bigg[\frac{1}{R^2}\, \Xi(R)\Bigg]^{\frac{1}{2}}\, .
\end{eqnarray}
We observe that \reff{cacc 2} implies (again recalling \reff{theo ass 1})
\begin{eqnarray}\label{proof 3}
\lefteqn{\int_{B_R} D^2f(\nabla u) \big(\nabla \partial_\alpha u, \nabla \partial_\alpha u\big) \Gamma^{-\frac{\gamma}{2}}
\gameps^{-\frac{\gamma_Q}{2}} \dx}\nonumber\\
&\leq &\frac{c}{R^2} \int_{T_R} \Gamma^{-\frac{\gamma + \mub}{2}} \gameps^{-\frac{\gamma_Q}{2}} 
|\partial _\alpha u - Q_\alpha|^2 \dx \, ,
\end{eqnarray}
hence we can make use of our assumption \reff{theo ass 2} by choosing a suitable subsequence $R \to \infty$ and obtain
\[
\int_{\rz^2} D^2f(\nabla u) \big(\nabla \partial_\alpha u, \nabla \partial_\alpha u\big) \Gamma^{-\frac{\gamma}{2}}
\gameps^{-\frac{\gamma_Q}{2}} \dx < \infty \, ,
\]
thus $I(R) \to 0$ as $R \to \infty$.\\

With this information we return to \reff{proof 1}, \reff{proof 2} and obtain
\[
\int_{B_r} \big| \nabla u - (\nabla u)\big|^2 \dx \leq I(R) \cdot \Bigg[\frac{1}{R^2}\, \Xi(R)\Bigg]^{\frac{1}{2}}
\]
with $I(R) \to 0$ as $R \to \infty$.  This proves the theorem with the help of hypothesis \reff{theo ass 2}. \qed\\

%***************************************************************************************************
%***************************************************************************************************
\section{Applications to the linear growth case}\label{linear}
%***************************************************************************************************
%***************************************************************************************************

%***** LABEL linear .... ******************************************************************************

Throughout this section we replace Assumption \ref{ass main 1} by a suitable stronger variant specifying the linear
growth condition. More precisely, we require

\begin{ass}\label{ass main linear}
The convex energy density $f$: $\rz^2 \to \rz$ is of class $C^2(\rz^2)$ and satisfies the non-uniform
ellipticity condition \reff{intro 2} with exponents $\mu > 1$, $\mub \leq 1$.

Moreover we assume that there exists a constant $M> 0$ such that for all $Z \in \rz^2$
\begin{equation}\label{linear 1}
|\nabla f(Z)| \leq M \, .
\end{equation}
\end{ass}

As outlined in \cite{BF:2021_1} (compare the discussion after inequality (1.3) in this reference),
Assumption \ref{ass main linear} implies with constants
$a$, $A >0$, $b$, $B \geq 0$ and for all $Z\in \rz^2$ the linear growth condition
\begin{equation}\label{linear 1a}
a |Z| - b \leq f(Z) \leq A |Z| + B\, . 
\end{equation} 

Before summarizing some corollaries of Theorem \ref{theo 1} in this particular setting,
we will show the following proposition which follows the line of the proof of Theorem 1.1 of \cite{BF:2021_1}.

\begin{proposition}\label{linear prop 1}
Suppose that we have Assumption \ref{ass main linear}. Then
\begin{equation}\label{linear 2}
\limsup_{|x| \to \infty} \frac{|u(x)|}{|x|} < \infty \, 
\end{equation}
implies
\begin{equation}\label{linear 3}
\int_{T_R} \Gamma^{\frac{1}{2}} \dx \leq c \big[1+R^2\big]\, .
\end{equation}
\end{proposition}
\emph{Proof of Proposition \ref{linear prop 1}.} By the convexity of $f$ we have for all $Z \in \rz^2$
\begin{equation}\label{linear 4}
f(Z)  \leq f(0) + \nabla f(Z)\cdot Z  \, ,
\end{equation}
hence for any $\eta \in C^1_0(\rz^2)$, $\eta \equiv 1$ on
$T_R$, $\op{spt} \eta \subset \hat{T}_R$, $0 \leq \eta \leq 1$, $|\nabla \eta| \leq c/R$,
we obtain using \reff{linear 1}, \reff{linear 1a} and \reff{linear 4}
\begin{eqnarray}\label{linear 5}
\int_{T_R} \Gamma^{\frac{1}{2}} \dx &\leq & c \int_{\hat{T}_R} \big[ 1+ f(\nabla u) \big]\eta^2 \dx\nonumber\\
& \leq & c  R^2 + c \int_{\hat{T}_R}  \nabla f(\nabla u)\cdot \nabla u  \eta^2 \dx  \, .
\end{eqnarray}
Now we use the weak form of equation \reff{intro 1} with testfunction $\psi = u \eta^2$, i.e.
\begin{eqnarray*}
0&=& \int_{\rz^2} \nabla f(\nabla u) \cdot \nabla \big[u \eta^2\big]\dx \nonumber\\
&=& \int_{\hat{T}_R} \nabla f(\nabla u) \cdot\nabla u \eta^2 \dx  + \int_{\hat{T}_R} \nabla f(\nabla u)\cdot \nabla \eta 2 \eta u \dx\, ,
\end{eqnarray*}
hence we have
\begin{equation}\label{linear  6}
\int_{\rz^2} \nabla f(\nabla u) \cdot \nabla u \eta^2\dx
\leq \frac{c}{R} \sup_{\hat{T}_R} |u| R^2 \, .
\end{equation}

Recalling our assumption \reff{linear 2} and combing \reff{linear 5} and \reff{linear 6} we find 
the claim \reff{linear 3} of the proposition.\qed\\

The first corollary to Theorem \ref{theo 1} immediately yields an extension of Theorem \ref{intro theo 1}. 

\begin{cor}\label{linear cor 1}
Theorem \ref{intro theo 1} remains valid for $1/2 < \mub \leq 1$.
\end{cor}

\emph{Proof of Corollary \ref{linear cor 1}.} We choose $\gamma_Q = 0$, $\gamma$ sufficiently close to $-1/2$ and let $Q =0$, i.e.
\[
\Xi(R) = \int_{T_R} \Gamma^{-\frac{1}{2}} |\nabla u|^2 \dx \leq \int_{T_R} |\nabla u|\dx \, ,
\]
and apply Proposition \ref{linear prop 1} to obtain hypothesis \reff{theo ass 2} of Theorem \ref{theo 1}. \qed\\

The next corollary gives an refinement of Corollary \ref{linear cor 1} by taking a measure for the relative oscillation
into account (compare \reff{linear cor_2 1} and \reff{linear cor_2 2}).

%***** LABEL linear cor_2 .... ******************************************************************************

\begin{cor}\label{linear cor 2} Suppose that we have Assumption \ref{ass main linear} with $1/2 < \mub < 1$.
For given $\gamma$, $\gamma_Q > 0$, $1-\mub < \gamma + \gamma_Q < 1/2$, we let
\[
p = \frac{1}{2-\gamma_Q - \gamma - \mub} > 1\, ,\quad q = \frac{1}{\gamma_Q+\gamma + \mub -1} 
\]
and define
\begin{equation}\label{linear cor_2 1}
\Theta_{Q}(x)  := \Bigg[\frac{\gameps}{\Gamma}\Bigg]^{\frac{2- \gamma_Q}{2} q}\, , \quad
\Theta(R) :=\inf_Q \frac{1}{|T_R|} \int_{T_R}\Theta_Q(x)\dx \leq 1 \, . 
\end{equation}

If we suppose that for all $R$ sufficiently large
\begin{equation}\label{linear cor_2 2}
\sup_{\hat{T}_R} |u|  \leq c R \Theta (R)^{-\frac{p}{q}} 
\end{equation}
with a constant $c > 0$ not depending on $R$, then $u$ is an affine function.
\end{cor}
\emph{Proof of Corollary \ref{linear cor 2}.} We estimate
\begin{eqnarray*}
\int_{T_R} \Gamma^{- \frac{\gamma+\mub}{2}} \gameps^{-\frac{\gamma_Q}{2}}|\nabla u -Q|^{2} \dx&\leq & 
\int_{T_R} \Theta_{Q}^{\frac{1}{q}} \Gamma^{\frac{2- \gamma_Q - \gamma -\mub}{2}}\dx 
\end{eqnarray*}
and obtain
\begin{eqnarray}\label{linear cor_2 3}
\Xi(R) &\leq & \int_{T_R} \Theta_{Q}^{\frac{1}{q}} \Gamma^{\frac{2- \gamma_Q - \gamma -\mub}{2}}\dx \nonumber\\
& \leq &
 \Bigg[\int_{T_R} \Theta_Q\dx\Bigg]^{\frac{1}{q}} \, \Bigg[\int_{T_R} \Gamma^{\frac{1}{2}}\dx\Bigg]^{\frac{1}{p}} \, .
\end{eqnarray}
We choose $Q = Q(R)$ such that 
\[
\frac{1}{|T_R|} \int_{T_R} \Theta_Q \dx \leq  2 \Theta(R) \, .
\]
Then \reff{linear cor_2 3} implies 
\begin{equation}\label{linear cor_2 4}
\Xi(R) \leq c \big[\Theta(R) R^2\big]^{\frac{1}{q}} \Bigg[\int_{T_R} \Gamma^{\frac{1}{2}}\dx \Bigg]^{\frac{1}{p}}\, .
\end{equation}
Discussing the right-hand side of \reff{linear cor_2 4} we exactly follow the proof of Proposition \ref{linear prop 1}
and just insert hypothesis \reff{linear cor_2 2} in \reff{linear 6}. This shows
\[
\Xi(R) \leq c \big[\Theta(R) R^2\big]^{\frac{1}{q}} \cdot \big[1+ \Theta(R)^{-\frac{p}{q}} R^2\big]^{\frac{1}{p}} \, ,
\]
hence we obtain Corollary \ref{linear cor 2} by recalling $\frac{1}{p}+\frac{1}{q}=1$. \qed\\

%***** LABEL linear cor_3 .... ******************************************************************************
A rather important class of energy densities with linear growth is of splitting-type, i.e.~of the form
\[
f(Z) = f_1(Z_1) + f_2 (Z_2)
\]
with functions $f_1$, $f_2$ of linear growth.  The particular features of splitting type energy densities with linear growth are discussed 
in \cite{BF:2020_3} (compare also \cite{BF:2020_4}).
If $f_1$, $f_2$ satisfy \reff{intro 2} with exponents $\mub_1$, $\mub_2 \leq 0$, respectively,
then we have $\mub = 0$ for the energy density $f$.\\

Nevertheless we still can derive Liouville-type theorems from Theorem \ref{theo 1}.  
In Corollary \ref{linear cor 3} we present an application, where in addition we make explicit use of
the flexibility of the vector $Q\in \rz^2$ by choosing $Q$ as a mean value. 
Then a smallness condition imposed on $|\nabla^2 u|$ provides the vanishing of the second derivatives.

\begin{cor}\label{linear cor 3} Suppose that we are given Assumption \ref{ass main linear} with
$\mub > -1/2$. If we have for a finite constant $c$ that
\begin{equation}\label{linear cor_3 1}
\sup_{T_R} |\nabla^2 u| \leq c R^{-1} \, ,
\end{equation}
then $u$ is an affine function.
\end{cor}
\emph{Proof of Corollary \reff{linear cor 3}.} We choose $\gamma = -\mub < 1/2$, $\gamma_Q =0$ such that \reff{theo ass 1} is
satisfied. We have to show that \reff{linear cor_3 1} implies \reff{theo ass 2}, i.e.~we claim that in this case
\begin{equation}\label{linear cor_3 2}
\Xi(R) = \int_{T_R} |\nabla u -Q|^2 \dx < c\big[1+R^2\big] \, ,
\end{equation}
where we choose $Q = (\nabla u)_R$. In fact, we have by the Poincar\'{e} inequality
(see \cite{St:1990_1}, Theorem A.10 ,  as the appropriate variant)
\begin{eqnarray*}
\int_{T_R} |\nabla u -Q|^2 \dx & \leq & c R^2 \int_{T_R} |\nabla^2 u|^2\dx
\leq  c  \Big[\sup_{T_R} |\nabla^2 u| \Big]^2 R^{4} \, ,
\end{eqnarray*}
which proves the corollary on account of the hypothesis  \reff{linear cor_3 1}. \qed\\

%***************************************************************************************************
%***************************************************************************************************
\section{Applications to the superlinear case}\label{pq}
%***************************************************************************************************
%***************************************************************************************************

We adapt Assumption \ref{ass main 1} to the case of superlinear growth, i.e.~we now require
\begin{ass}\label{ass main pq}
Suppose that we are given numbers $\mu\in \rz$, $q>1$ such that $-\mu \leq q-2$ and let $\mub := 2-q$ be the exponent
on the right-hand side of \reff{intro 2}, i.e.~the convex energy density $f$: $\rz^2 \to \rz$ is of class $C^2(\rz^2)$ and satisfies the non-uniform ellipticity condition
\begin{equation}\label{pq 1}
c_1 \big(1+|Z|^2\big)^{-\frac{\mu}{2}} |Y|^2 \leq D^2 f(Z)\big(Y,Y\big) \leq c_2 \big(1+|Z|^2\big)^{\frac{q-2}{2}}|Y|^2
\end{equation}
for all $Z$, $Y\in \rz^2$ with exponents $\mu > 1$, $q > 1$ and constants $c_1$, $c_2 >0$.

Suppose that we have in addition
\begin{equation}\label{pq 2}
a |Z|^s - b \leq f(Z) \quad\mbox{for some $1 < s \leq q$}
\end{equation}
and with constants $a>0$, $b \geq 0$.
\end{ass}

Note that, as in the linear growth case, the auxiliary parameter $\mu$ needs no further specification in our hypotheses.\\

Conditions \reff{pq 1} and \reff{pq 2} are introduced in \cite{BF:2001_1} describing energy densities of $(s,\mu,q)$-growth, we refer
to \cite{Bi:1818}, Section 3.2, for a more detailed discussion. In particular \reff{pq 1} implies with some constant $M>0$ and for all $Z \in \rz^2$
\begin{equation}\label{pq 3}
|\nabla f(Z) | \leq M \big(1+|Z|^{2}\big)^{\frac{q-1}{2}} \, .
\end{equation} 
Examples are given, for instance, by
\begin{eqnarray*}
f_1(Z) &=& \big(1+|Z_1|^2\big)^{\frac{s}{2}} + |Z_2|^2 \, \quad 1 < s \leq 2\, ,\\[2ex]
f_2(Z) &=& \big(1+|Z|^2\big)^{\frac{s}{2}} + \big(1+|Z_2|\big)^{\frac{q}{2}}\, \quad 1 < s \leq q\, . 
\end{eqnarray*}

The main result of this section is 
\begin{cor}\label{pq cor 1}
Given Assumption \ref{ass main pq} we suppose in addition that
\begin{equation}\label{pq 4}
s > q - \frac{1}{2}\, .
\end{equation}
If we have
\begin{equation}\label{pq 5}
\limsup_{|x|\to \infty}\frac{|u(x)|}{|x|} < \infty \, ,
\end{equation}
then $u$ is an affine function.
\end{cor}

\emph{Proof of Corollary \ref{pq cor 1}.}
In order to apply Theorem \ref{theo 1} we let $Q=0$, $\gamma_Q =0$ and since we have \reff{pq 4} we can choose 
$0 < \gamma < 1/2$ such that
\begin{equation}\label{pq 6}
s >  q- \gamma \, .
\end{equation}
Then, in our main Theorem \ref{theo 1} we observe
\begin{equation}\label{pq 7}
\Xi(R) \leq  \int_{\hat{T}_R} \Gamma^{\frac{q-\gamma}{2}} \dx \, .
\end{equation} 
We follow the proof of Proposition \ref{linear prop 1}, where now \reff{linear 5} and \reff{linear 6} are replaced by (recalling \reff{pq 2}, \reff{pq 3} 
\reff{pq 5} and choosing $l \in \nz$ sufficiently large)
\begin{eqnarray}\label{pq 8}
\int_{\hat{T}_R}  \Gamma^{\frac{s}{2}} \eta^{2l}\dx &\leq& c \int_{\hat{T}_R} \big[1+ f(\nabla u)\big] \eta^{2l}\dx\nonumber\\[2ex]
&\leq & c R^2 + \frac{c}{R} \sup_{\hat{T}_R} |u|\int_{\hat{T}_R} |\nabla f(\nabla u)| \eta^{2l-1} \dx \nonumber\\[2ex]
& \leq & c R^2 + c \int_{\hat{T}_R} \Gamma^{\frac{q-1}{2}} \eta^{2l-1}\dx \, .
\end{eqnarray}
Since $q-1 < q-\gamma$ we find real numbers $\hat{q}$, $\hat{p} > 1$, $\frac{1}{\hat{q}} + \frac{1}{\hat{p}} =1$ such that
\begin{equation}\label{pq 9}
\int_{\hat{T}_R} \Gamma^{\frac{q-1}{2}}\eta^{2l-1} \dx \leq c
\Bigg[ \int_{\hat{T}_R} \Gamma^{\frac{q-\gamma}{2}}\eta^{2l}\dx\Bigg]^{\frac{1}{\hat{q}}} R^{\frac{2}{\hat{p}}}\, .
\end{equation}
From \reff{pq 6}, \reff{pq 8} and \reff{pq 9} we obtain
\begin{eqnarray}\label{pq 10}
\int_{\hat{T}_R} \Gamma^{\frac{q-\gamma}{2}}\eta^{2l}\dx &\leq &c  \int_{\hat{T}_R} \Gamma^{\frac{s}{2}} \eta^{2l} \dx\nonumber\\
&\leq & c R^2 + c \Bigg[\int_{\hat{T}_R} \Gamma^{\frac {q-\gamma}{2}} \eta^{2l}\dx\Bigg]^{\frac{1}{\hat{q}}} R^{\frac{2}{\hat{p}}}\, .
\end{eqnarray}
W.l.o.g.~we suppose that
\[
R^2 \leq c \Bigg[\int_{\hat{T}_R} \Gamma^{\frac {q-\gamma}{2}} \eta^{2l}\dx\Bigg]^{\frac{1}{\hat{q}}} R^{\frac{2}{\hat{p}}}
\]
and \reff{pq 10} yields
\begin{equation}\label{pq 11}
\Bigg[\int_{\hat{T}_R} \Gamma^{\frac{q-\gamma}{2}} \eta^{2l} \dx\Bigg]^{1-\frac{1}{\hat{q}}} \leq c R^{\frac{2}{\hat{p}}}\, .
\end{equation}
Recalling \reff{pq 7} and $1-\frac{1}{\hat{q}} = \frac{1}{\hat{p}}$ we have proved \reff{theo ass 2} and this finally implies Corollary \ref{pq cor 1}. \qed\\

\bibliography{Cacc_Liou}
\bibliographystyle{unsrt}

\end{document}